\documentclass{article} % For LaTeX2e
\usepackage{nips15submit_e,times}

\usepackage{hyperref}
\usepackage{graphicx} % to include figur\usepackage{kasai}es
\usepackage{subfigure}
\usepackage{enumitem}
\usepackage{epstopdf}

\usepackage{amssymb}
\usepackage{amsmath}
\usepackage{amsthm}
\usepackage{mathtools}
\usepackage[raggedright]{titlesec} % avoid hyphen in sections

\usepackage{algorithm,algcompatible,amsmath}
\algnewcommand\INPUT{\item[\textbf{Input:}]}%
\algnewcommand\OUTPUT{\item[\textbf{Output:}]}%

\usepackage{arydshln}
\newcommand{\rank}{\mathrm{rank}}

\newcommand{\Stiefel}[2]{{\mathrm{St}({#1},{#2})}}
\newcommand{\Grass}[2]{{\mathrm{Gr}({#1},{#2})}}

\newcommand{\OG}[1]{{\mathcal{O}({#1})}}

\newcommand{\grad}{\mathrm{grad}}

\newcommand{\mat}[1]{{\bf #1}}

\newcommand{\Grad}{\mathrm{Grad}}
\newcommand{\Exp}{{\mathrm{Exp}}}
\newcommand{\Log}{{\mathrm{Log}}}
\renewcommand{\cos}{{\mathrm{cos}}}
\renewcommand{\sin}{{\mathrm{sin}}}

\makeatletter
\def\hlinewd#1{%
  \noalign{\ifnum0=`}\fi\hrule \@height #1 \futurelet
   \reserved@a\@xhline}
\makeatother

\title{A Riemannian gossip approach to\\ decentralized matrix completion}

\date{Compiled on \today}

\author{
Bamdev Mishra\\
Amazon Development Centre India,\\
Bangalore, India \\
\texttt{bamdevm@amazon.com}  \\
\And
Hiroyuki Kasai\\
The University of Electro-Communications\\
Tokyo, Japan\\
\texttt{kasai@is.uec.ac.jp} \\
\And
Atul Saroop\\
Amazon Development Centre India,\\
Bangalore, India\\
\texttt{asaroop@amazon.com}\\
}

\nipsfinalcopy % Uncomment for camera-ready version

\begin{document}

%%  BM: 8-page manuscript + references

\maketitle

\begin{abstract}
In this paper, we propose novel gossip algorithms for the low-rank decentralized matrix completion problem. The proposed approach is on the Riemannian Grassmann manifold that allows \emph{local} matrix completion by different agents while achieving asymptotic consensus on the \emph{global} low-rank factors. The resulting approach is scalable and parallelizable. Our numerical experiments show the good performance of the proposed algorithms on various benchmarks.
\end{abstract}

\section{Introduction}
The problem of low-rank matrix completion amounts to completing a matrix from a small number of entries by assuming a low-rank model for the matrix. The problem has many applications in control systems and system identification \cite{markovsky13a}, collaborative filtering \cite{rennie05a}, and information theory \cite{shi16a}, to name a just few. Consequently, it has been a topic of great interest and there exist many large-scale implementations for both \emph{batch} \cite{cai10a, toh10a, ngo12a, mishra14c, wen12a, boumal15a, keshavan10a} and \emph{online} scenarios that focus on parallel and stochastic implementations \cite{balzano10a, yu14a, recht13a, teflioudi12a}.

In this paper, we are interested in a \emph{decentralized} setting, where we divide the matrix completion problem into smaller subproblems that are solved by many agents locally while simultaneously enabling them to arrive at a \emph{consensus} that solves the full problem \cite{lin15a}. The recent paper \cite{lin15a} proposes a particular decentralized framework for matrix completion by exploiting the algorithm proposed in \cite{ngo12a}. It, however, requires an inexact dynamic consensus step at every iteration. We relax this by proposing a novel formulation that combines together a weighted sum of completion and consensus terms. Additionally, in order to minimize the communication overhead between the agents, we constrain each agent to communicate with \emph{only} one other agent as in the \emph{gossip} framework \cite{boyd06a}. One motivation is that this addresses privacy concerns of sharing sensitive data \cite{lin15a}. Another motivation is that the gossip framework is robust to scenarios where certain agents may be inactive at certain time slots, e.g., consider each agent to be a computing machine. We propose a preconditioned variant that is particularly well suited for ill-conditioned instances. Additionally, we also propose a parallel variant that allows to exploit parallel computational architectures. All the variants come with asymptotic convergence guarantees. To the best of our knowledge, this is the first work that exploits the gossip architecture for solving the decentralized matrix completion problem.

The organization of the paper is as follows. In Section \ref{sec:problem_setup}, we discuss the decentralized problem setup and propose a novel problem formulation. In Section \ref{sec:gossip}, we discuss the proposed stochastic gradient gossip algorithm for the matrix completion problem. A preconditioned variant of the Riemannian gossip algorithm is motivated in Section \ref{sec:preconditioned_variant}. Additionally, we discuss a way to parallelize the proposed algorithms in Section \ref{sec:parallel_variant}. Numerical comparisons in Section \ref{sec:comparisons} show that the proposed algorithms compete effectively with state-of-the-art on various benchmarks. The Matlab codes for the proposed algorithms are available at \url{https://bamdevmishra.com/codes/gossipMC/}.

\section{Decentralized matrix completion}\label{sec:problem_setup}
The matrix completion problem is formulated as
\begin{equation}\label{eq:matrix-completion-formulation}
\begin{array}{llll}
	\min\limits_{\mat{X}\in\mathbb{R}^{m \times n}}
		& \displaystyle\frac{1}{2}	\|\mathcal{P}_{\Omega}(\mat{X}) - \mathcal{P}_{\Omega}(\mat{X}^\star)\|_F^2 \\
	 \text{subject to} & \rank(\mat{X})=r,
\end{array}
\end{equation}
where ${\mat{X }^ \star}\in\mathbb{R}^{n\times m}$ is a matrix whose entries are known for indices if they belong to the subset $(i,j)\in\Omega$ and $\Omega$ is a subset of the complete set of indices $\{(i,j):i\in\{1,...,m\}\text{ and }j\in\{1,...,n\}\}$. The operator $\mathcal{P}_{\Omega}(\mat{X}_{ij})=\mat{X}_{ij}$ if $(i,j) \in \Omega$ and $\mathcal{P}_{\Omega}(\mat{X}_{ij})=0$ otherwise is called the orthogonal sampling operator and is a mathematically convenient way to represent the subset of known entries. The rank constraint parameter $r$ is usually set to a low value, i.e., $\ll (m, n)$ that implies that we seek low-rank completion. A way to handle the rank constraint in (\ref{eq:matrix-completion-formulation}) is by a fixed-rank matrix parameterization. In particular, we use $\mat{X} = \mat{U} \mat{W}^T$, where $\mat{U} \in \Stiefel{r}{m}$ and $\mat{W} \in \mathbb{R}^{n \times r}$, where $ \Stiefel{r}{m}$ is the set of $m\times r$ \emph{orthonormal} matrices, i.e., the columns are orthonormal. The interpretation is that $\mat{U}$ captures the dominant column space of $\mat X$ and $\mat W$ captures the \emph{weights} \cite{mishra14a}. Consequently, the optimization problem (\ref{eq:matrix-completion-formulation}) reads 
\begin{equation}\label{eq:factorized_formulation}
\min_{\mat{U}\in \Stiefel{r}{m}}  \min_{\mat{W}\in\mathbb{R}^{n \times r }} \|\mathcal{P}_{\Omega}(\mat{UW}^T) - \mathcal{P}_{\Omega}(\mat{X^\star})\|_F^2.
\end{equation}
The \emph{inner} least-squares optimization problem in (\ref{eq:factorized_formulation}) is solved in closed form by exploiting the least-squares structure of the cost function to obtain the optimization problem
\begin{equation}\label{eq:grassmann_formulation}
\begin{array}{lll}
\min\limits_{\mat{U}\in\Stiefel{r}{m}}  & \displaystyle\frac{1}{2}\|\mathcal{P}_{\Omega}(\mat{UW_{\mat U}}^T) - \mathcal{P}_{\Omega}({ \mat X}^\star)\|_F^2
\end{array}
\end{equation}
in $\mat U$, where $\mat{W}_{\mat U}$ is the solution to the inner optimization problem $\min_{\mat{W}\in\mathbb{R}^{n \times r }}  \|\mathcal{P}_{\Omega}(\mat{UW}^T) - \mathcal{P}_{\Omega}(\mat{X}^\star)\|_F^2$. (The cost function in (\ref{eq:grassmann_formulation}) may be discontinuous at points $\mat U$ where $\mat{W}_{\mat U}$ is non-unique \cite{dai12a}. This is handled effectively by adding a regularization term $\| \mat{X} \|_F^2$ to (\ref{eq:matrix-completion-formulation}).)
%and $f: \Stiefel{r}{m} \rightarrow \mathbb{R}$. 

The problem (\ref{eq:grassmann_formulation}) requires handling the entire incomplete matrix $\mat{X}^\star$ at all steps of optimization. This is memory intensive and computationally heavy, especially in large-scale instances. To relax this constraint, we distribute the task of solving the problem (\ref{eq:grassmann_formulation}) among $N$ agents, which perform certain computations independently. To this end, we partition the incomplete matrix $\mat{X}^\star = [\mat{X}_1^\star, \mat{X}_2^\star,\ldots, \mat{X}_N^\star]$ along the columns such that the size of $\mat{X}_i^\star$ is $m\times n_i$ with $\sum n_i = n$ for $i=\{1,2,\ldots, N\}$. Each agent $i$ has knowledge of the incomplete matrix $\mat{X}_i^\star$ and its \emph{local} set of indices $\Omega_i$ of known entries. We also partition the weight matrix $\mat W$ as $\mat{W}^T = [\mat{W}_1^T, \mat{W}_2^T, \ldots, \mat{W}_N ^T]$ such that the matrix $\mat{W}_i$ has size $n_i\times r$. A straightforward reformulation of (\ref{eq:grassmann_formulation}) is 
\begin{equation}\label{eq:semidecentralized_grassmann_formulation}
\begin{array}{lll}
 \displaystyle \sum\limits_{i} \min\limits_{\mat{U}\in\Stiefel{r}{m}, \mat{W}_i \in  \mathbb{R}^{n_i \times r}}  &   \displaystyle\frac{1}{2} \|\mathcal{P}_{\Omega _i}({{\mat{UW}}_i}^T) - \mathcal{P}_{\Omega _i}(\mat{X}^\star _i)\|_F^2 \\
= \min\limits_{\mat{U}\in\Stiefel{r}{m}}  &  \displaystyle\frac{1}{2} \sum\limits_{i} \underbrace{\|\mathcal{P}_{\Omega _i}({{\mat{UW}}}_{i{\mat U}}^T) - \mathcal{P}_{\Omega _i}(\mat{X}^\star _i)\|_F^2}_{{\rm problem\ handled\ by\ agent\ }i},
\end{array}
\end{equation}
where ${{\mat{W}}}_{i{\mat U}}$ is the least-squares solution to $\min_{\mat{W}_i \in\mathbb{R}^{n_i \times r }}  \|\mathcal{P}_{\Omega _i}(\mat{UW}_i^T) - \mathcal{P}_{\Omega _i}(\mat{X}_i^\star)\|_F^2$, which can be computed by agent $i$ independently of other agents. 

Although the computational workload gets distributed among the agents in the problem formulation (\ref{eq:semidecentralized_grassmann_formulation}), all agents require the knowledge of $\mat U$ (to compute matrices $ {\mat{W}}_{i{\mat U}}$). To circumvent this issue, instead of one shared matrix $\mat U$ for all agents, each agent $i$ stores a local copy $\mat{U}_i$, which it then updates based on information from its \emph{neighbors}. For minimizing the communication overhead between agents, we additionally put the constraint that at any time slot only two agents communicate, i.e, each agent has exactly only one neighbor. This is the basis of the gossip framework \cite{boyd06a}. In standard gossip framework, at a time slot, an agent is randomly assigned one neighbor \cite{boyd06a}. However, to motivate the various ideas in this paper and to keep the exposition simple, we fix the agents network topology, i.e., each agent is preassigned a unique neighbor. (In Section \ref{sec:arbritrary_topology}, we show how to deal with random assignments of neighbors.) To this end, the agents are numbered according to their proximity, e.g., for $ i \leqslant N -1$, agents $i$ and $i + 1$ are neighbors. Equivalently, agents $1$ and $2$ are neighbors and can communicate. Similarly, agents $2$ and $3$ communicate, and so on. This communication between the agents allows to reach a \emph{consensus} on $\mat{U}_i$. Specifically, it suffices that the \emph{column spaces} of all $\mat{U}_i$ converge. (The precise motivation and formulation are in Section \ref{sec:gossip}.) Our proposed decentralized matrix completion problem formulation is 
\begin{equation}\label{eq:decentralized_grassmann_formulation}
\begin{array}{lll}
 \min\limits_{\mat{U}_1 , \ldots, \mat{U}_N \in\Stiefel{r}{m}}  & \displaystyle\frac{1}{2} \sum\limits_{i} \underbrace{\|\mathcal{P}_{\Omega _i}({{\mat{U}_i\mat{W}}}_{i{\mat{U}_i}}^T) - \mathcal{P}_{\Omega _i}(\mat{X}^\star _i)\|_F^2}_{{\rm completion\ task\ handled\ by\ agent\ } i} \\ 
 & \qquad +\  \displaystyle\frac{\rho}{2} \underbrace{( d_1 (\mat{U}_1, \mat{U}_2)^2 + d_2 (\mat{U}_2, \mat{U}_3)^2 + \ldots + d_{N-1} (\mat{U}_{N-1}, \mat{U}_{N})^2 )}_{\rm consensus},
\end{array}
\end{equation}
where $d_i$ is a certain distance measure (defined in Section \ref{sec:gossip}) between $\mat{U}_i$ and $\mat{U}_{i+1}$ for $i\leqslant N-1$, minimizing which forces $\mat{U}_i$ and $\mat{U}_{i+1}$ to an ``average'' point (specifically, average of the column spaces). $\rho \geqslant 0$ is a parameter that trades off matrix completion with consensus. Here $\mat{W}_{i{\mat U}_i}$ is the solution to the optimization problem $\min_{\mat{W}_i \in \mathbb{R}^{n_i \times r }}  \|\mathcal{P}_{\Omega _i }(\mat{U}_i \mat{W}_i^T) - \mathcal{P}_{\Omega _i}(\mat{X}_i ^\star)\|_F^2$. 

In standard gossip framework, the aim is to make the agents converge to a common point, e.g, minimizing only the consensus term in (\ref{eq:decentralized_grassmann_formulation}). In our case, we additionally need the agents to perform certain tasks, e.g., minimizing the completion term in (\ref{eq:decentralized_grassmann_formulation}), which motivates the \emph{weighted} formulation (\ref{eq:decentralized_grassmann_formulation}). For a large value of $\rho$, the consensus term in (\ref{eq:decentralized_grassmann_formulation}) dominates, minimizing which allows the agents to arrive at consensus. For $\rho = 0$, the optimization problem (\ref{eq:decentralized_grassmann_formulation}) solves $N$ independent completion problems and there is no consensus. For a sufficiently large value of $\rho$, the problem (\ref{eq:decentralized_grassmann_formulation}) achieves the goal of matrix completion along with consensus.

%The Riemannian distance $d$ captures the distance between two subspaces. 
%
%
%This allows the agents to come to a consensus
%
%The basis for this is the \emph{gossip framework} on the Grassmann manifold $\Grass{r}{m}$ \cite[Section~4.4]{ bonnabel13a}.

\section{The Riemannian gossip algorithm}\label{sec:gossip}
It should be noted that the optimization problem (\ref{eq:grassmann_formulation}), and similarly (\ref{eq:semidecentralized_grassmann_formulation}), only depends on the \emph{column space} of $\mat U$ rather than $\mat U$ itself \cite{boumal15a, balzano10a}. Equivalently, the cost function in (\ref{eq:grassmann_formulation}) remains constant under the transformation $\mat U \mapsto \mat{UO}$ for all orthogonal matrices $\mat O$ of size $r \times r$. Mathematically, the column space of $\mat U$ is captured by the set, called the \emph{equivalence class}, of matrices 
\begin{equation}\label{eq:eqivalence_class}
[\mat U]:=\{ \mat{UO}: \mat{O} {\rm \ is\ a\ }r\times r {\rm \ orthogonal\ matrix}\}.
\end{equation}
The set of the equivalence classes is called the \emph{Grassmann} manifold, denoted by $\Grass{r}{m}$, which is the set of $r$-dimensional subspaces in $\mathbb{R}^m$ \cite{absil08a}. The Grassmann manifold $\Grass{r}{m}$ is identified with the quotient manifold $\Stiefel{r}{m}/\OG{r}$, where $\OG{r}$ is the orthogonal group of $r\times r$ matrices \cite{absil08a}. 

Subsequently, the problem (\ref{eq:grassmann_formulation}), and similarly (\ref{eq:semidecentralized_grassmann_formulation}), is on the Grassmann manifold $\Grass{r}{m}$ and not on $\Stiefel{r}{m}$. However, as $ \Grass{r}{m}$ is an abstract quotient space, numerical optimization algorithms are implemented with matrices $\mat U$ on $\Stiefel{r}{m}$, but conceptually, optimization is on $\Grass{r}{m}$. It should be stated that the Grassmann manifold has the structure of a \emph{Riemannian} manifold and optimization on the Grassmann manifold is a well studied topic in literature. Notions such as the Riemannian gradient (first order derivatives of a cost function), geodesic (shortest distance between elements), and logarithm mapping (capturing ``difference'' between elements) have closed-form expressions \cite{absil08a}.

\begin{table}[t]
\caption{Proposed online gossip algorithm for (\ref{eq:final_decentralized_grassmann_formulation})}
\label{tab:online_algorithm} 
\begin{center} \small
\begin{tabular}{ |p{13.5cm}| }
\hline
\begin{enumerate}
\item At each time slot $t$, pick an agent $i \leqslant N-1$ randomly with uniform probability.
\item Compute the Riemannian gradients $\grad_{x_i} f_i$, $\grad_{x_{i+1}} f_{i + 1}$, $\grad_{x_i} d_i$, and $\grad_{x_{i+1}} d_{i}$ with the matrix representations
\[
\begin{array}{lll}
\Grad_{x_i}f_i =(\mathcal{P}_{\Omega _i}({{\mat{U}_i\mat{W}}}_{i{\mat{U}_i}}^T) - \mathcal{P}_{\Omega _i}(\mat{X}^\star _i) )\mat{W}_{i{\mat{U}_i}} \\
\grad_{x_i}f_i = \Grad_{x_i}f_i - \mat{U}_i (\mat{U}_i ^T \Grad_{x_i}f_i) \\
\grad_{x_i} d_i = -\Log_{x_i}(x_{i + 1})\\
\grad_{x_{i+1}} d_i = -\Log_{x_{i+1}}(x_{i}),
\end{array}
\]
where $\mat{U}_i$ is the matrix representation of $x_i$. $\Log_{x_i}(x_{i + 1})$ is the \emph{logarithm} mapping, which is defined as
\[
\begin{array}{lll}
\Log_{x_i}(x_{i + 1}) = \mat{P} {\rm arctan}(\mat{S})\mat{Q}^T,
\end{array}
\]
where $\mat{PS Q}^T$ is the rank-$r$ singular value decomposition of $( \mat{U}_{i+1} -   \mat{U}_i  (\mat{U}_{i}^T \mat{U}_{i+1})) \break   (\mat{U}_{i}^T \mat{U}_{i+1})^{-1} $. The $\rm{arctan}(\cdot)$ operation is only on the diagonal entries. It should be noted that the Riemannian  gradient of the Riemannian distance is the negative logarithm mapping \cite{huper10a}.
\item Given a stepsize $\gamma_t$, update $x_i$ and $x_{i + 1}$ as 
\[
\begin{array}{lllll}
{x_i}_+  = \Exp_{x_i} (- \gamma_t (\alpha_i \grad_{x_i} f_i + \rho    \grad_{x_i} d_i ) ) \\
{x_{i + 1}}_+  =  \Exp_{x_{i+1}}( - \gamma_t (  \alpha_{i+1} \grad_{x_{i+1}} f_{i + 1} + \rho \grad_{x_{i+1}} d_{i})), \\
\end{array}
\]
where $\mat{U}_i$ is the matrix representation of $x_i$ and $\alpha_i = 1$ if $i=\{1, N \}$, else $\alpha_i = 0.5$. $\Exp_{x_i}(\xi_{x_i}) = \mat{U}_{i} \mat{V} \cos(\Sigma)\mat{V}^T +   \mat{W} \sin(\Sigma) \mat{V}^T$ is the \emph{exponential} mapping and $\mat{W \Sigma V}^T$ is the rank-$r$ singular value decomposition of $\xi_{x_i}$. The $\cos(\cdot)$ and $\sin(\cdot)$ operations are only on the diagonal entries.
%\item Repeat.
\end{enumerate}
  \\
  \\  
 \hline
\end{tabular}
\end{center} 
\end{table}

If $x$ is an element of a Riemannian compact manifold $\mathcal{M}$, then the decentralized formulation (\ref{eq:decentralized_grassmann_formulation}) boils down to the form
\begin{equation}\label{eq:final_decentralized_grassmann_formulation}
\begin{array}{lll}
\min\limits_{x_1, \ldots, x_N \in \mathcal{M}}  & \displaystyle \sum_i {f_i(x_i)} \\ 
   & \qquad +\  \displaystyle\frac{\rho}{2} \underbrace{( d_1 (x_1, x_2)^2 + d_2 (x_2, x_3)^2 + \ldots + d_{N-1} (x_{N-1}, x_{N})^2 )}_{\rm consensus},
\end{array}
\end{equation}
where $x_i =[\mat{U} _i]$ with matrix representation ${\mat U}_i \in \Stiefel{r}{m}$, $\mathcal{M} = \Grass{r}{m} = \Stiefel{r}{m}/\OG{r}$, $f_i: \mathcal{M} \rightarrow \mathbb{R}$ is a continuous function, and $d_i: \mathcal{M} \times \mathcal{M} \rightarrow \mathbb{R}$ is the Riemannian geodesic distance between $x_i$ and $x_{i+1}$. Here $[\mat{U}_i]$ is the equivalence class defined in (\ref{eq:eqivalence_class}). The Riemannian distance $d_i$ captures the distance between the subspaces $[\mat{U}_i]$ and $[\mat{U}_{i + 1}]$. Minimizing only the consensus term in (\ref{eq:final_decentralized_grassmann_formulation}) is equivalent to computing the \emph{Karcher mean} of $N$ subspaces \cite{bonnabel13a, huper10a}.

We exploit the stochastic gradient descent setting framework proposed by Bonnabel \cite{bonnabel13a} for solving (\ref{eq:final_decentralized_grassmann_formulation}), which is an optimization problem on the Grassmann manifold. In particular, we exploit the stochastic gradient algorithm in the gossip framework. To keep the analysis simple, we predefine the topology on the agents network. Following \cite[Section~4.4]{ bonnabel13a}, we make the following assumptions.
\begin{enumerate}[label=\textbf{A\arabic*}]
\item \label{item:A1} Agents $i$ and $i + 1$ are neighbors for all $i \leqslant N-1$.

\item  \label{item:A2} At each time slot, say $t$, we pick an agent $i \leqslant N-1$ randomly with uniform probability. This means that we also pick agent $i+1$ (the neighbor of agent $i$). Subsequently, agents $i$ and $i + 1$ update $x_i$ and $x_{i +1}$, respectively, by taking a \emph{gradient descent step} with stepsize $\gamma_t$ on $\mathcal{M}$. The stepsize sequence satisfies the standard conditions, i.e., $\sum\limits \gamma_t ^2 < \infty$ and $\sum\limits\gamma_t = +\infty$ \cite[Section~3]{bonnabel13a}. 
\end{enumerate}

Each time we pick an agent $i \leqslant N-1$, we equivalently also pick its neighbor $i+1$. Subsequently, we need to update both of them by taking a gradient descent step based on $f_i(x_i) + f_{i+1}(x_{i+1}) + \rho d_i(x_i, x_{i+1})^2 /2$. Repeatedly updating the agents in this fashion is a stochastic process.

It should be noted that because of the particular topology and sampling that we assume (in \ref{item:A1} and \ref{item:A2}), on an average $x_2$ to $x_{N-1}$ are updated \emph{twice} the number of times $x_1$ and $x_N$ are updated. For example, if $N=3$, then \ref{item:A1} and \ref{item:A2} lead to solving (in expectation) the problem $\min_{x_1, x_2, x_3 \in \mathcal{M}}  f_1(x_1) + 2f_2(x_2)+f_3(x_3) + \rho (d_1(x_1, x_2)^2 + d_2(x_2, x_3)^2)/2$. To resolve this issue, we multiply the scalar $\alpha_i$ to $f_i$ (and its Riemannian gradient) while updating $x_i$s. Specifically, $\alpha_i = 1$ if $i=\{1, N \}$, else $\alpha_i = 0.5$. If $\grad_{x_i} f_i$ is the Riemannian gradient of $f_i$ at $x_i \in \mathcal{M}$, then the stochastic gradient descent algorithm updates $x_i$ along the search direction $-(\alpha_i \grad_{x_i} f_i + \rho \grad_{x_i} d_i)$ with the exponential mapping $\Exp_{x_i} : T_{x_i} \mathcal{M} \rightarrow \mathcal{M}$, where $T_{x_i}\mathcal{M}$ is the tangent space of $\mathcal{M}$ at $x_i$. The overall algorithm with concrete matrix expressions is in Table \ref{tab:online_algorithm}. The stochastic gradient descent algorithm in Table \ref{tab:online_algorithm} converges to a critical point of (\ref{eq:final_decentralized_grassmann_formulation}) \emph{almost surely} \cite{bonnabel13a}. The gradient updates require the computation of the Riemannian gradient of the cost function in (\ref{eq:final_decentralized_grassmann_formulation}) and moving along the geodesics with \emph{exponential} mapping, both of which have closed-form expressions on the Grassmann manifold $\Grass{r}{m}$ \cite{absil08a}. Similarly, the matrix completion problem specific gradient computations are shown in \cite{boumal15a}.

\subsection{Computational complexity}
For an update of $x_i$ with the formulas shown in Table \ref{tab:online_algorithm}, the computational complexity depends on the computation of \emph{partial derivatives} of the cost function in (\ref{eq:final_decentralized_grassmann_formulation}), e.g., $\Grad_{x_i} f_i$. Particularly, in the context of the problem (\ref{eq:decentralized_grassmann_formulation}), the computational cost is $O(|\Omega_i|r^2 + n_i r ^2 + m r)$. The Grassmann manifold related ingredients, e.g., $\Exp$, cost $O(mr^2 + r^3)$.

\subsection{Convergence analysis}\label{sec:convergence_analysis}
Asymptotic convergence analysis of the algorithm in Table \ref{tab:online_algorithm} follows directly from the analysis in \cite[Theorem~1]{bonnabel13a}. The key idea is that for a compact Riemannian manifold all continuous functions of the parameter can be bounded. This is the case for (\ref{eq:final_decentralized_grassmann_formulation}), which is on the \emph{compact} Grassmann manifold $\Grass{r}{m}$. Subsequently, under a decreasing stepsize condition and noisy gradient estimates (that is an unbiased estimator of the batch gradient), the stochastic gradient descent algorithm in Table \ref{tab:online_algorithm} converges to a critical point of (\ref{eq:final_decentralized_grassmann_formulation}) almost surely. Conceptually, while the standard stochastic gradient descent setup deals with an \emph{infinite} stream of samples, we deal with a finite number of samples (i.e., we pick an agent $i \leqslant N-1$), which we repeat many times.

\begin{table}[t]
\caption{Proposed preconditioned gossip algorithm for (\ref{eq:final_decentralized_grassmann_formulation})}
\label{tab:precon_online_algorithm} 
\begin{center} \small
\begin{tabular}{ |p{13.5cm}| }
\hline
\begin{enumerate}
\item At each time slot $t$, pick an agent $i \leqslant N-1$ randomly with uniform probability and compute the Riemannian gradients $\grad_{x_i} f_i$, $\grad_{x_{i+1}} f_{i + 1}$, $\grad_{x_i} d_i$, and $\grad_{x_{i+1}} d_{i}$ with the matrix representations shown in Table \ref{tab:online_algorithm}.
\item Given a stepsize $\gamma_t$, update $x_i$ and $x_{i + 1}$ as 
\[
\begin{array}{lllll}
{x_i}_+  = \Exp_{x_i} (- \gamma_t (\alpha_i \grad_{x_i} f_i + \rho    \grad_{x_i} d_i ) ({\mat{W}_{i{\mat U}_i} ^T \mat{W}_{i{\mat U}_i}} + {\rho \mat{I}})^{-1}) \\
{x_{i + 1}}_+  =  \Exp_{x_{i+1}}( - \gamma_t (  \alpha_{i+1} \grad_{x_{i+1}} f_{i + 1} + \rho \grad_{x_{i+1}} d_{i}) ({\mat{W}_{{i+1}{\mat U}_{i+1}} ^T \mat{W}_{{i+1}{\mat U}_{i+1}}} + {\rho \mat{I}})^{-1} ), \\
\end{array}
\]
where $\mat{W}_{i{\mat U}_i}$ is the least-squares solution to the optimization problem $\min_{\mat{W}_i \in \mathbb{R}^{n_i \times r }}  \|\mathcal{P}_{\Omega _i }(\mat{U}_i \mat{W}_i^T) - \mathcal{P}_{\Omega _i}(\mat{X}_i ^\star)\|_F^2$. $\Exp$ and $\alpha_i$ are defined in Table \ref{tab:online_algorithm}.
%\item Repeat.
\end{enumerate}
  \\
 % \\  
 \hline
\end{tabular}
\end{center} 
\end{table}

\subsection{Preconditioned variant}\label{sec:preconditioned_variant}
The performance of first order algorithm (including stochastic gradients) often depends on the \emph{condition number} (the ratio of maximum eigenvalue to the minimum eigenvalue) of the Hessian of the cost function (at the minimum). The issue of ill-conditioning arises especially when data $\mat{X}^\star$ have drawn power law distributed singular values. Additionally, a large value of $\rho$ in (\ref{eq:final_decentralized_grassmann_formulation}) leads to convergence issues for numerical algorithms. To this end, the recent works \cite{ngo12a, mishra14c, boumal15a} exploit the concept of \emph{manifold preconditioning} in matrix completion. Specifically, the Riemannian gradients are \emph{scaled} by computationally cheap matrix terms that arise from the second order curvature information of the cost function. Matrix scaling of the gradients is equivalent to multiplying an approximation of the \emph{inverse} Hessian to gradients. This operation on a manifold requires special attention. In particular, the matrix scaling \emph{must} be a positive definite operator on the tangent space of the manifold \cite{mishra14c, boumal15a}.

Given the Riemannian gradient $\xi_{x_i} = \grad_{x_i} f_i + \rho \grad_{x_i} d_i$ computed by agent $i$, the proposed manifold preconditioning is
\begin{equation}\label{eq:preconditioner}
\xi_{x_i} \mapsto \xi_{x_i}   (\underbrace{\mat{W}_{i{\mat U}_i} ^T \mat{W}_{i{\mat U}_i}}_{\rm from\ completion}  + \underbrace{\rho \mat{I}}_{\rm from\ consensus})^{-1},
\end{equation}
where $\mat{W}_{i{\mat U}_i}$ is the solution to the optimization problem $\min_{\mat{W}_i \in \mathbb{R}^{n_i \times r }}  \|\mathcal{P}_{\Omega _i }(\mat{U}_i \mat{W}_i^T) - \mathcal{P}_{\Omega _i}(\mat{X}_i ^\star)\|_F^2$ and $\mat{I}$ is $r\times r$ identity matrix. The use of preconditioning (\ref{eq:preconditioner}) costs $O(n_i r^2 + r^3)$.

The term $\mat{W}_{i{\mat U}_i} ^T \mat{W}_{i{\mat U}_i}$ is motivated by the fact that it is computationally cheap to compute and captures a \emph{block diagonal approximation} of the Hessian of the simplified (but related) cost function $ \|{{\mat{U}_i\mat{W}}}_{i{\mat{U}_i}}^T - \mat{X}^\star _i \|_F^2$. The works \cite{ngo12a, mishra14c, boumal15a} use such \emph{preconditioners} with superior performance. The term $\rho \mat{I}$ is motivated by the fact that the second order derivative of the square of the Riemannian geodesic distance is an identity matrix. Finally, it should be noted that the matrix scaling is positive definite, i.e., ${\mat{W}_{i{\mat U}_i} ^T \mat{W}_{i{\mat U}_i}}  + {\rho \mat{I}} \succ 0$ and that the transformation (\ref{eq:preconditioner}) is on the tangent space. Equivalently, if $\xi_{x_i}$ belongs to $T_{x_i} \mathcal{M}$, then $\xi_{x_i} (\mat{W}_{i{\mat U}_i} ^T \mat{W}_{i{\mat U}_i}  + \rho \mat{I})^{-1}$ also belongs to $T_{x_i} \mathcal{M}$. This is readily checked by the fact that the tangent space $T_{x_i} \mathcal{M}$ at $x_i$ on the Grassmann manifold is characterized by the set $\{ \eta_{x_i}:\eta_{x_i} \in \mathbb{R}^{m \times r},\mat{U}_{i}^T \eta_{x_i} = 0 \}$.

The proposed preconditioned variant of the stochastic gradient descent algorithm for (\ref{eq:final_decentralized_grassmann_formulation}) is shown in Table \ref{tab:precon_online_algorithm}. It should be noted that preconditioning the gradients does not affect the asymptotic convergence guarantees of the proposed algorithm.

\subsection{Parallel variant}\label{sec:parallel_variant}
Assumption \ref{item:A1} on the network topology of agents allows to propose parallel variants of the proposed stochastic gradient descent algorithms in Tables \ref{tab:online_algorithm} and \ref{tab:precon_online_algorithm}. To this end, instead of picking one agent at a time, we pick agents in such a way that it leads to a number of parallel updates. 

We explain the idea for $N = 5$. Updates of the agents are divided into two \emph{rounds}. In round $1$, we pick agents $1$ and $3$, i.e., all the \emph{odd} numbered agents. It should be noted that the neighbor of agent $1$ is agent $2$ and the neighbor of agent $3$ is agent $4$. Consequently, the updates of agents $1$ and $2$ are independent from those of agents $3$ and $4$ and hence, can be carried out in parallel. In round $2$, we pick agents $2$ and $4$, i.e., all the \emph{even} numbered agents. The updates of agents $2$ and $3$ are independent from those of agents $4$ and $5$ and therefore, can be carried out in parallel.

\begin{table}[t]
\caption{Proposed parallel variant for (\ref{eq:final_decentralized_grassmann_formulation})}
\label{tab:parallel_algorithm} 
\begin{center} \small
\begin{tabular}{ |p{13.5cm}| }
\hline
\begin{enumerate}
\item Define round $1$ as consisting of agents $i = {1, 3, \ldots}$ and their neighbors. Define round $2$ as consisting of agents $i = {2, 4, \ldots}$ and their neighbors.
\item At each time slot $t$, pick a round $j \leqslant 2$ randomly with uniform probability.
\item Given a stepsize, update the agents (and their corresponding neighbors) in parallel with the updates proposed in Table \ref{tab:online_algorithm} (or in Table \ref{tab:precon_online_algorithm}).
%\item Repeat.
\end{enumerate}
  \\
 % \\  
 \hline
\end{tabular}
\end{center} 
\end{table}

The key idea is that \emph{randomness} is on the rounds and not on the agents. For example, we pick a round $j$ from $\{1, 2 \}$ with uniform probability. Once a round is picked, the updates on the agents (that are part of this round) are performed with the same stepsize and in parallel. The stepsize is updated when a new round is picked. The stepsize sequence satisfies the standard conditions, i.e., it is square-summable and its summation is divergent. The overall algorithm in shown in Table \ref{tab:parallel_algorithm}.

To prove convergence, we define two new functions,
\begin{equation}\label{eq:new_functions}
\begin{array}{lll}
g_1(x_1, x_2, \ldots) = f_1(x_1) + f_2(x_2)+\ldots + \displaystyle\frac{\rho}{2} ( d_1 (x_1, x_2)^2 + d_3 (x_3, x_4)^2 +  \ldots)\\
g_2(x_2, x_3, \ldots) = f_2(x_2) + f_3(x_3)+\ldots + \displaystyle\frac{\rho}{2} ( d_2 (x_2, x_3)^2 + d_4 (x_4, x_5)^2 +  \ldots),
\end{array}
\end{equation}
that consist of terms from the cost function in (\ref{eq:final_decentralized_grassmann_formulation}). The algorithm in Table \ref{tab:parallel_algorithm} is then interpreted as the standard stochastic gradient descent algorithm applied to the problem 
\begin{equation} \label{eq:parallel_formulation}
\min_{x_i \in \mathcal{M}} \quad g_1(x_1, x_2, \ldots) + g_2(x_2, x_3, \ldots).
\end{equation}
with two ``samples'' that are chosen randomly at each time slot. Consequently, following the standard arguments, the algorithm in Table \ref{tab:parallel_algorithm} converges asymptotically to a critical point of (\ref{eq:parallel_formulation}). However, it should also be noted that the addition of $g_1$ and $g_2$ leads to $x_2$ to $x_{N-1}$ being updated (on an average) \emph{twice} the number of times $x_1$ and $x_N$ are updated. This is handled by multiplying $\alpha_i$ to $f_i$ while updating $x_i$s, where $\alpha_i = 1$ if $i=\{1, N \}$, else $\alpha_i = 0.5$. Finally, the algorithm in Table (\ref{tab:parallel_algorithm}) converges to a critical point of (\ref{eq:final_decentralized_grassmann_formulation}). It should emphasized that parallelization of the updates is for free by virtue of construction of functions in (\ref{eq:new_functions}).

\subsection{Extension to continuously changing network topology}\label{sec:arbritrary_topology}
The algorithm in Table \ref{tab:online_algorithm} assumes that the neighbors of the agents are predefined in a particular way (assumption \ref{item:A1}). However, in many scenarios the network topology changes with time \cite{boyd06a}. To simulate the scenario, we first consider a fully connected network of $N$ agents. The number of unique edges is $N(N-1)/2$. We pick an edge $ik$ (the edge that connects agents $i$ and $k$) randomly with uniform probability and drop all the other edges. Equivalently, only one edge is \emph{active} at any time slot. Consequently, we update agents $i$ and $k$ with a gradient descent update, e.g., based on Table \ref{tab:online_algorithm} or Table \ref{tab:precon_online_algorithm}. The overall algorithm is shown in Table \ref{tab:arbitrary_online_algorithm}. Following the arguments in Section \ref{sec:convergence_analysis}, it is straightforward to see that the proposed algorithm converges almost surely to a critical point of a problem that combines completion along with consensus, i.e., 
\begin{equation}\label{eq:arbitrary_decentralized_grassmann_formulation}
\begin{array}{lll}
\min\limits_{x_1, \ldots, x_N \in \mathcal{M}}  & \displaystyle (N-1)\sum_i f_i(x_i) +\  \displaystyle\frac{\rho}{2} { \sum_{i < k} d_{i  k} (x_i, x_k)^2 },
\end{array}
\end{equation}
where $d_{ik}(x_i, x_k)$ is the Riemannian geodesic distance between $x_i$ and $x_k$.

\begin{table}[t]
\caption{Proposed algorithm for continuously changing network topology}
\label{tab:arbitrary_online_algorithm} 
\begin{center} \small
\begin{tabular}{ |p{13.5cm}| }
\hline 
\begin{enumerate}
\item At each time slot $t$, pick a pair of agents, say $i$ and $k$, randomly with uniform probability.
\item Compute the Riemannian gradients $\grad_{x_i} f_i$, $\grad_{x_{k}} f_{k}$, $\grad_{x_i} d_{ik}$, and $\grad_{x_{k}} d_{ik}$ as
\[
\begin{array}{lll}
\Grad_{x_i}f_i =(\mathcal{P}_{\Omega _i}({{\mat{U}_i\mat{W}}}_{i{\mat{U}_i}}^T) - \mathcal{P}_{\Omega _i}(\mat{X}^\star _i) )\mat{W}_{i{\mat{U}_i}} \\
\grad_{x_i}f_i = \Grad_{x_i}f_i - \mat{U}_i (\mat{U}_i ^T \Grad_{x_i}f_i) \\
\grad_{x_i} d_{ik} = -\Log_{x_i}(x_{k}) \\
 \grad_{x_{k}} d_{ik} = -\Log_{x_{k}}(x_{i}),
\end{array}
\]
where $\Log$ is defined in Table \ref{tab:online_algorithm}.
\item Given a stepsize $\gamma_t$, update $x_i$ and $x_{k}$ as 
\[
\begin{array}{lllll}
{x_i}_+  = \Exp_{x_i} (- \gamma_t ( \grad_{x_i} f_i + \rho    \grad_{x_i} d_{ik} ) ) \\
{x_{k}}_+  =  \Exp_{x_{k}}( - \gamma_t(  \grad_{x_{k}} f_{k} + \rho \grad_{x_{k}} d_{ik})), \\
\end{array}
\]
where the exponential mapping $\Exp_{x_i}$ is defined in Table \ref{tab:online_algorithm}.
%\item Repeat.
\end{enumerate}\\
 \hline
\end{tabular}
\end{center} 
\end{table}

\section{Numerical comparisons}\label{sec:comparisons}
Our proposed algorithms in Table \ref{tab:online_algorithm} (Online Gossip) and in Table \ref{tab:precon_online_algorithm} (Precon Online Gossip) and their parallel variants, Parallel Gossip and Precon Parallel Gossip, are compared on different problem instances. The implementations are based on the Manopt toolbox \cite{boumal14a} with certain operations relying on the mex files supplied with \cite{boumal15a}. We also show comparisons with D-LMaFit, the decentralized algorithm proposed in \cite{lin15a} on smaller instances as the D-LMaFit code (supplied by the authors) is \emph{not} tuned to large-scale instances. As the mentioned algorithms are well suited for different scenarios, we compare them against the number of \emph{updates} performed by the agents. We fix the number of agents $N$ to $6$. Online algorithms are run for a maximum of $1000$ iterations. The parallel variants are run for $400$ iterations. Overall, agents $1$ and $N$ perform a maximum of $200$ updates (rest all perform $400$ updates). D-LMaFit is run for $400$ iterations, i.e., each agent performs $400$ updates. Algorithms are initialized randomly. The stepsize sequence is defined as $\gamma_t = \gamma_0/t$, where $t$ is the time slot and  $\gamma_0$ is set using cross validation. For simplicity, all figures only show the plots for agents $1$ and $2$.

All simulations are performed in Matlab on a $2.7$ GHz Intel Core i$5$ machine with $8$ GB of RAM. For each example considered here, an $m \times n$ random matrix of rank $r$ is generated as in \cite{cai10a}. Two matrices $\mat{A} \in \mathbb{R}^{m \times r}$ and $\mat{B} \in \mathbb{R}^{n \times r}$ are generated according to a Gaussian distribution with zero mean and unit standard deviation. The matrix product $\mat{AB} ^T$ gives a random matrix of rank $r$. A fraction of the entries are randomly removed with uniform probability and noise (sampled from the Gaussian distribution with mean zero and standard deviation $10^{-6} $) is added to each entry to construct the training set $\Omega$ and $\mat{X}^\star$. The over-sampling ratio (OS) is the ratio of the number of known entries to the matrix dimension, i.e, ${\rm OS} = |\Omega|/(mr +nr -r^2)$. We also create a test set by randomly picking a small set of entries from $\mat{AB} ^T$. The matrices $\mat{X}^\star _i$ are created by distributing the number of $n$ columns of $\mat{X}^\star$ equally among the agents. The training and test sets are also partitioned similarly.

\begin{figure*}[h]
\begin{tabular}{cc}
\hspace*{-0.8cm}
\noindent \begin{minipage}[b]{0.35\hsize}
\centering
\includegraphics[width=\hsize]{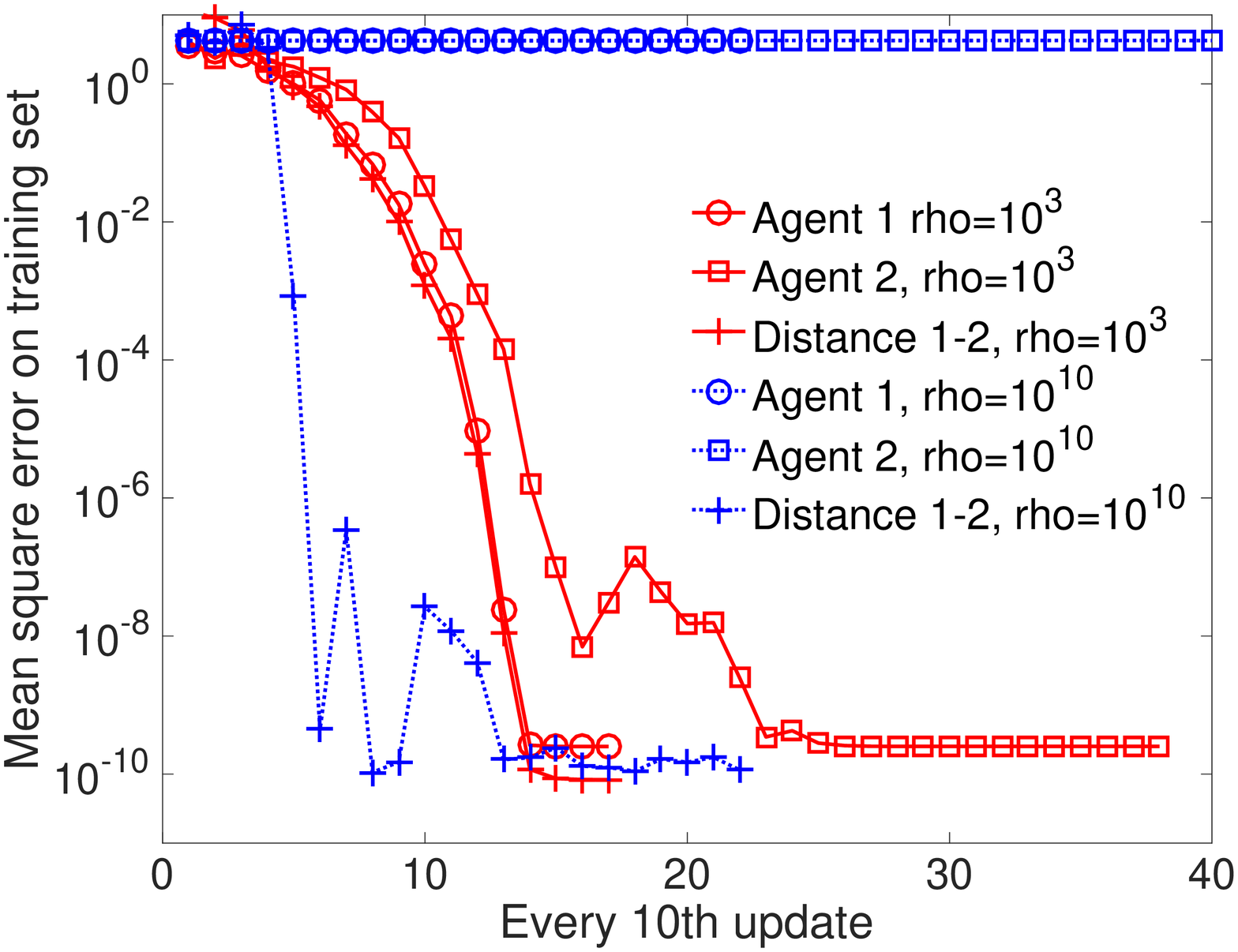}
{\scriptsize (a) Effect of $\rho$.}
\end{minipage}
%\hspace*{0.5cm}
\noindent \begin{minipage}[b]{0.35\hsize}
\centering
\includegraphics[width=\hsize]{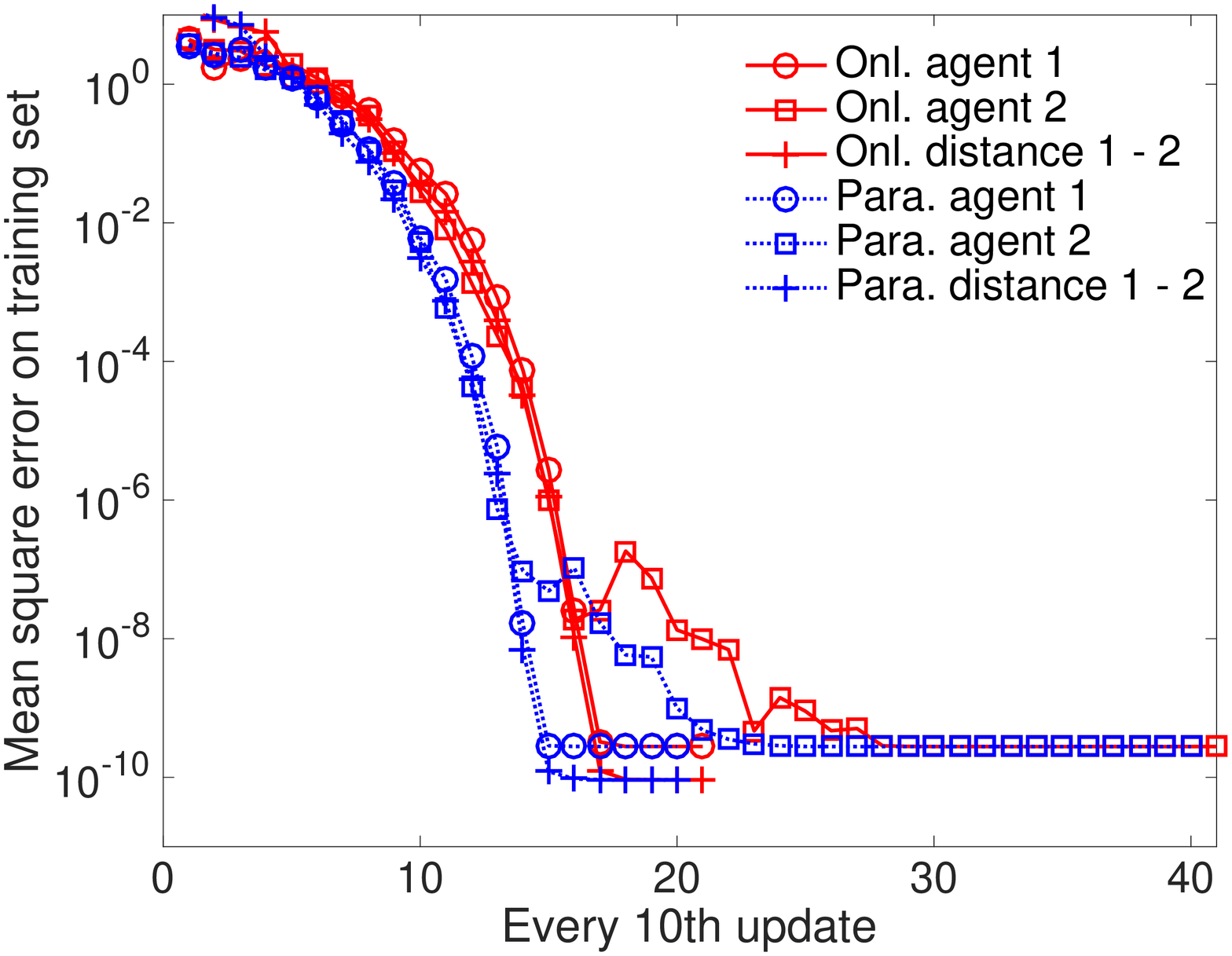}
{\scriptsize (b) Performance of online and parallel variants.}
 \end{minipage}
%\hspace*{-0.5cm}
\noindent \begin{minipage}[b]{0.35\hsize}
\centering
\includegraphics[width=\hsize]{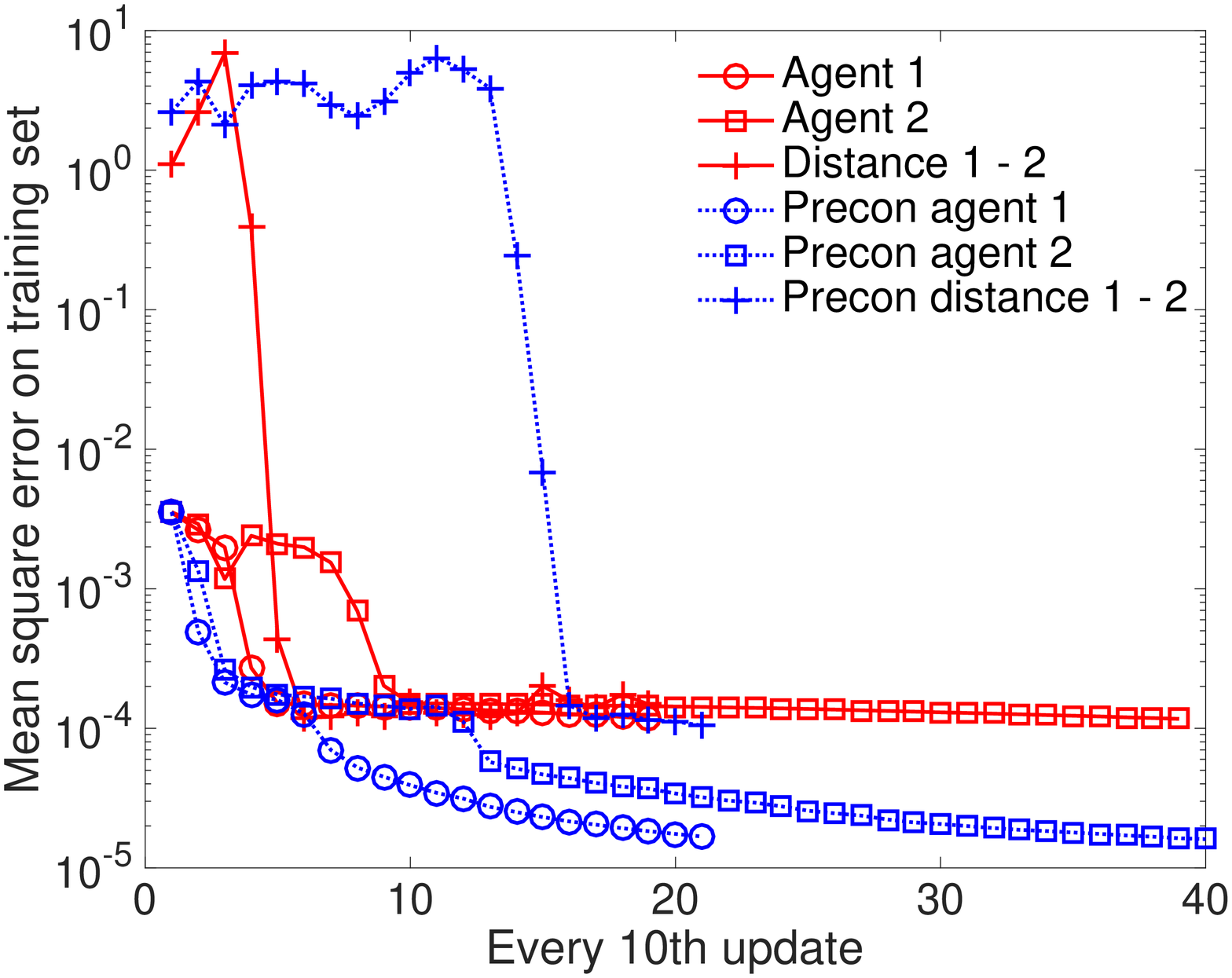}
{\scriptsize (c) Effect of preconditioning on training.}
\end{minipage}\\
\hspace*{-0.8cm}
\noindent \begin{minipage}[b]{0.35\hsize}
\centering
\includegraphics[width=\hsize]{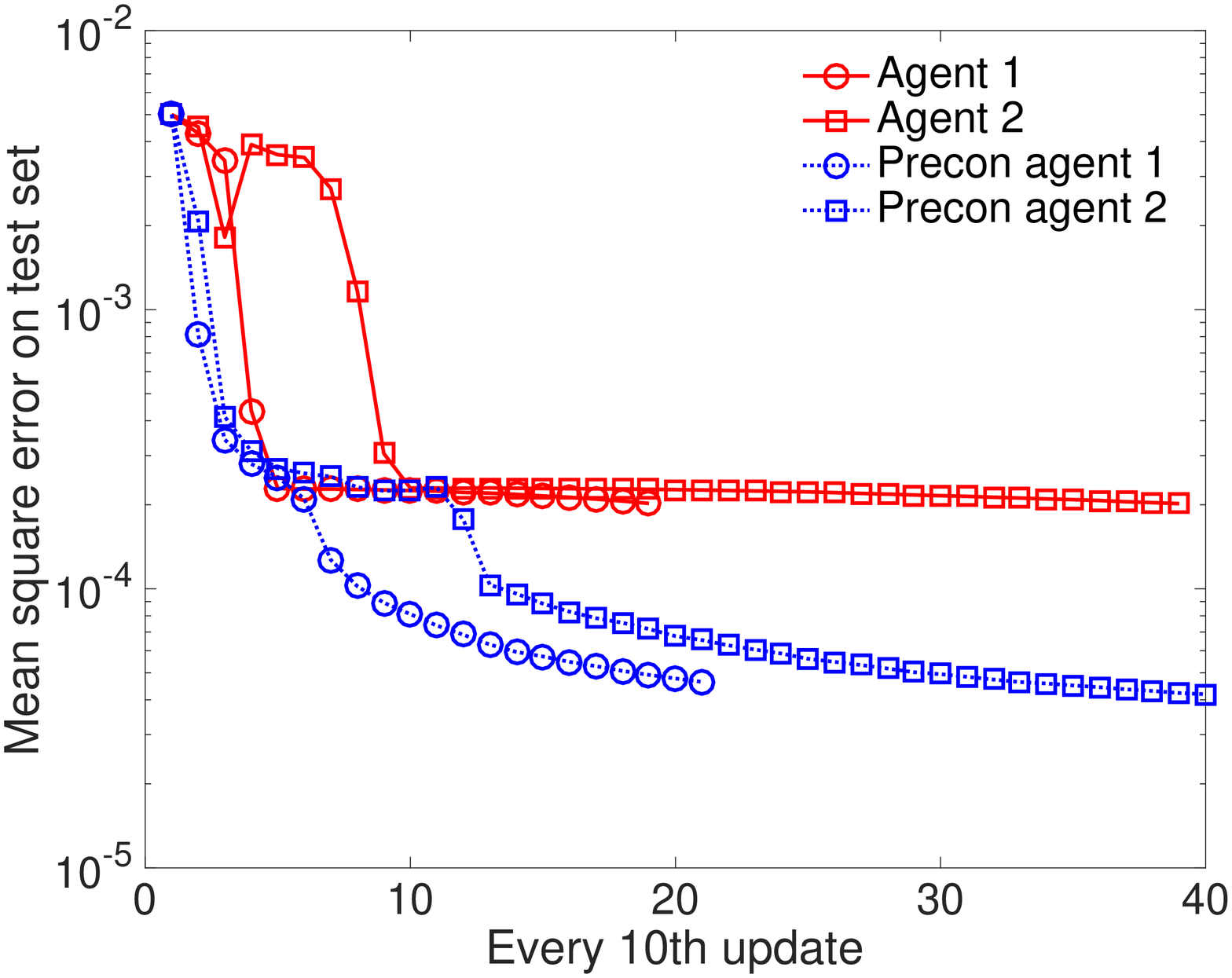}
{\scriptsize (d) Effect of preconditioning on test error.}
\end{minipage}
% \hspace*{-0.5cm}
\noindent \begin{minipage}[b]{0.35\hsize}
\centering
\includegraphics[width=\hsize]{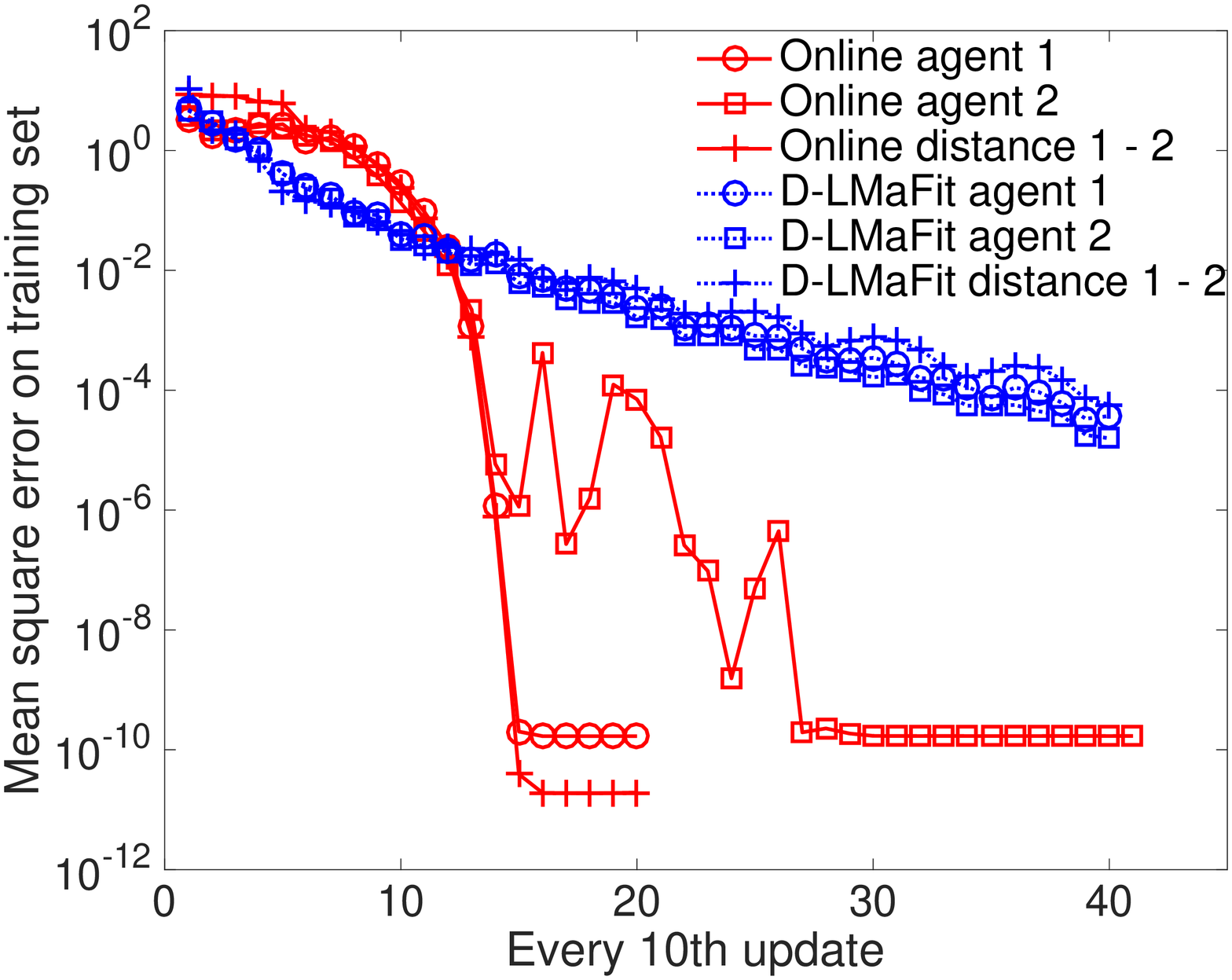}
{\scriptsize (e) Online Gossip outperforms D-LMaFit.}
\end{minipage}
\noindent \begin{minipage}[b]{0.35\hsize}
\centering
\includegraphics[width=\hsize]{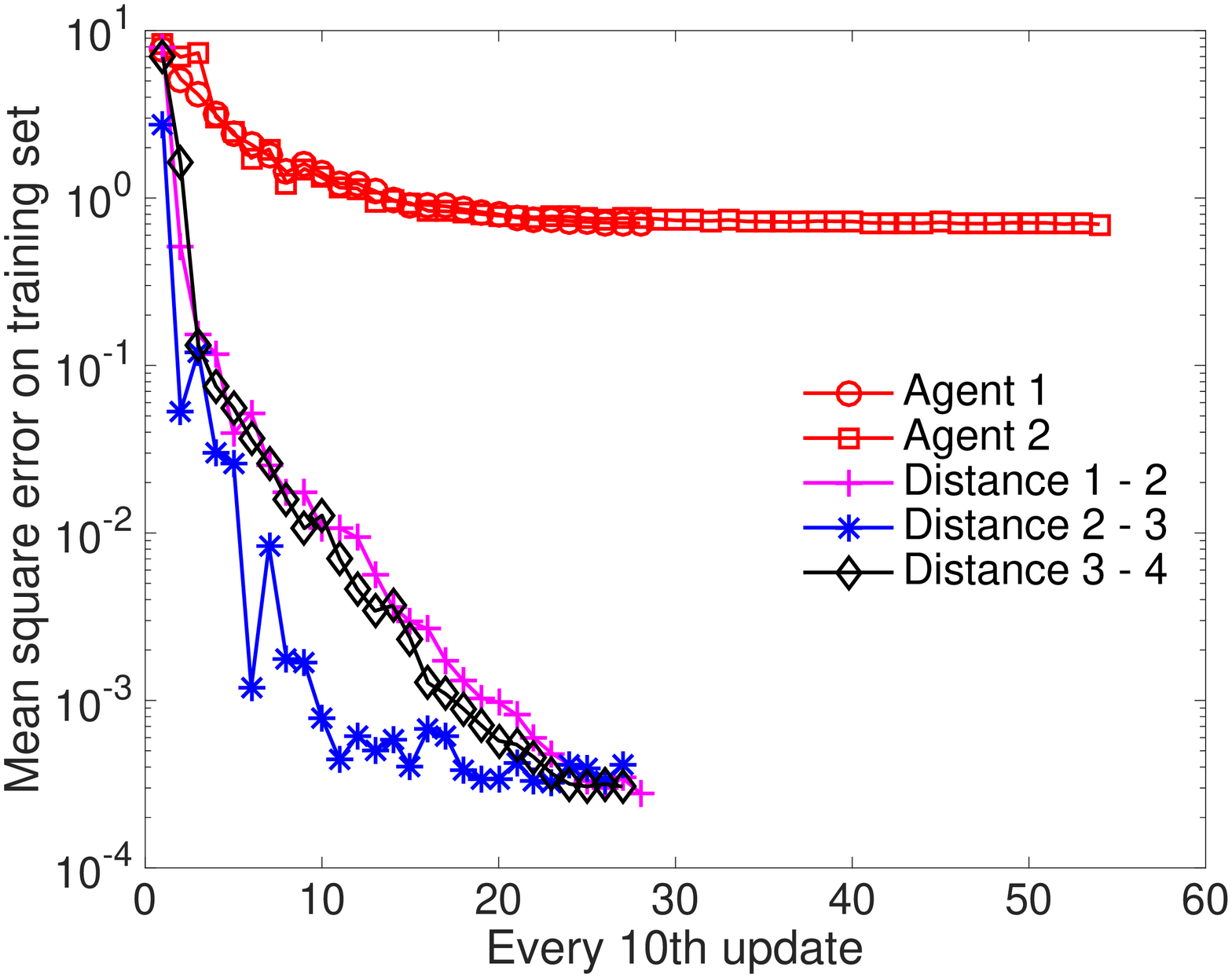}
{\scriptsize (f) MovieLens 20 M: consensus of the agents.}
\end{minipage}
\end{tabular}
\caption{Performance of proposed algorithms in different scenarios.}
\label{fig:all}
\end{figure*}

\textbf{Case 1: effect of $\rho$.} We consider a problem instance of size 10 000$\times $100 000 of rank $5$ and OS $6$. Two scenarios with $\rho = 10^3$ and $\rho = 10^{10}$ are considered. Figure \ref{fig:all}(a) shows the performance of Online Gossip. Not surprisingly, for $\rho = 10^{10}$, we only see consensus (the distance between agents $1$ and $2$ tends to zero). For $\rho = 10^3$, we see both completion and consensus, which validates the theory. 

\textbf{Case 2: performance of online versus parallel.} We consider Case 1 with $\rho = 10^3$. Figure \ref{fig:all}(b) shows the performance of Online Gossip and Parallel Gossip, both of which show a similar behavior on the training and test (not shown here) sets. 

\textbf{Case 3: ill-conditioned instances.} We consider a problem instance of size 5 000$\times $50 000 of rank $5$ and impose an exponential decay of singular values with condition number $500$ and OS $6$. Figure \ref{fig:all}(c) shows the performance of Online Gossip and its preconditioned variant for $\rho = 10^3$. During the initial updates, the preconditioned variant aggressively minimizes the completion term of (\ref{eq:decentralized_grassmann_formulation}), which shows the effect of the preconditioner (\ref{eq:preconditioner}). Eventually, consensus among the agents is achieved. Overall, the preconditioned variant shows a superior performance on both the training and test sets as shown in Figures \ref{fig:all}(c) and \ref{fig:all}(d). 

\textbf{Case 4: Comparisons with D-LMaFit \cite{lin15a}.}
We consider a problem instance of size $500\times 12000$, rank $5$, and OS $6$. D-LMaFit is run with the default parameters. For Online Gossip, we set $\rho = 10^3$. As shown in Figure \ref{fig:all}(e), Online Gossip quickly outperforms D-LMaFit. Overall, Online Gossip takes fewer number of updates to reach a high accuracy.

\textbf{Case 5: MovieLens 20M dataset \cite{harper15a}.}
Finally, we show the performance of Online Gossip on the MovieLens-20M dataset of $20000263$ ratings by $138493$ users for $26744$ movies. (D-LMaFit is not compared as it does not scale to this dataset.) We perform $5$ random $80/20$ train/test partitions. We split both the train and test data among $N=4$ agents along the number of users such that each agent has ratings for $26744$ movies and $34624$ (except agent $4$, which has $34621$) unique users. This ensures that the ratings are distributed evenly among the agents. We run Online Gossip with $\rho = 10^7$ (through cross validation) and for $800$ iterations. Figure \ref{fig:all}(f) shows that asymptotic consensus is achieved among the four agents. It should be noted that the distance between agents $2$ and $3$ decreases faster than others as agents $2$ and $3$ are updated (on an average) twice the number of times than agents $1$ and $4$ (assumption \ref{item:A1}). Table \ref{tab:movielens} shows the normalized mean absolute errors (NMAE) obtained on the \emph{full} test set averaged over five runs. NMAE is defined as the mean absolute error (MAE) divided by variation of the ratings. Since the ratings vary from $0.5$ to $5$, NMAE is MAE/$4.5$. We obtain the lowest NMAE for rank $5$.

\begin{table}[H]
\caption{Performance of Online Gossip on MovieLens 20M dataset}
\label{tab:movielens} 
\begin{center} \small
\begin{tabular}{ p{2.3cm}|p{2.3cm}|p{2.5cm}|p{2.3cm}|p{2.3cm}}
\hline
& Rank $3$ & Rank $5$ & Rank $7$ & Rank $9$\\
\hdashline
NMAE on test set & $0.1519 \pm 3\cdot 10^{-3}$ & $ \bf{ 0.1507 \pm 3\cdot 10^{-3}}$ &$0.1531 \pm 2\cdot 10^{-3}$ & $0.1543 \pm 1\cdot 10^{-3}$\\
%$9$ & $0.1386$ \\
 \hline
\end{tabular}
\end{center} 
\end{table}

\section{Conclusion}
We have proposed a Riemannian gossip approach to the decentralized matrix completion problem. Specifically, the completion task is distributed among a number of agents, which are then required to achieve consensus. Exploiting the gossip framework, this is modeled as minimizing a weighted sum of \emph{completion} and \emph{consensus} terms on the Grassmann manifold. The rich geometry of the Grassmann manifold allowed to propose a novel stochastic gradient descent algorithm for the problem with simple updates. Additionally, we have proposed two variants -- preconditioned and parallel -- of the algorithm for dealing with different scenarios. Numerical experiments show the competitive performance of the proposed algorithms on different benchmarks.

%\clearpage

\small
\bibliographystyle{unsrt}  
\bibliography{gossipMC_MKS2016}

\begin{thebibliography}{10}

\bibitem{markovsky13a}
I.~Markovsky and K.~Usevich.
\newblock Structured low-rank approximation with missing data.
\newblock {\em SIAM Journal on Matrix Analysis and Applications},
  34(2):814--830, 2013.

\bibitem{rennie05a}
J.~Rennie and N.~Srebro.
\newblock Fast maximum margin matrix factorization for collaborative
  prediction.
\newblock In {\em International Conference on Machine learning (ICML)}, pages
  713--719, 2005.

\bibitem{shi16a}
J.~Shi, Y.and~Zhang and K.~B. {Letaief}.
\newblock Low-rank matrix completion for topological interference management by
  {R}iemannian pursuit.
\newblock {\em IEEE Transactions on Wireless Communications}, PP(99), 2016.

\bibitem{cai10a}
J.~F. Cai, E.~J. Cand\`es, and Z.~Shen.
\newblock A singular value thresholding algorithm for matrix completion.
\newblock {\em SIAM Journal on Optimization}, 20(4):1956--1982, 2010.

\bibitem{toh10a}
K.~C. Toh and S.~Yun.
\newblock An accelerated proximal gradient algorithm for nuclear norm
  regularized least squares problems.
\newblock {\em Pacific Journal of Optimization}, 6(3):615--640, 2010.

\bibitem{ngo12a}
T.~T. Ngo and Y.~Saad.
\newblock {Scaled gradients on Grassmann manifolds for matrix completion}.
\newblock In {\em Advances in Neural Information Processing Systems 25 (NIPS)},
  pages 1421--1429, 2012.

\bibitem{mishra14c}
B.~Mishra and R.~Sepulchre.
\newblock {R3MC}: A {R}iemannian three-factor algorithm for low-rank matrix
  completion.
\newblock In {\em Proceedings of the 53rd IEEE Conference on Decision and
  Control (CDC)}, pages 1137--1142, 2014.

\bibitem{wen12a}
Z.~Wen, W.~Yin, and Y.~Zhang.
\newblock Solving a low-rank factorization model for matrix completion by a
  nonlinear successive over-relaxation algorithm.
\newblock {\em Mathematical Programming Computation}, 4(4):333--361, 2012.

\bibitem{boumal15a}
N.~Boumal and P.-A. Absil.
\newblock Low-rank matrix completion via preconditioned optimization on the
  {G}rassmann manifold.
\newblock {\em Linear Algebra and its Applications}, 475:200--239, 2015.

\bibitem{keshavan10a}
R.~H. Keshavan, A.~Montanari, and S.~Oh.
\newblock Matrix completion from a few entries.
\newblock {\em IEEE Transactions on Information Theory}, 56(6):2980--2998,
  2010.

\bibitem{balzano10a}
L.~Balzano, R.~Nowak, and B.~Recht.
\newblock Online identification and tracking of subspaces from highly
  incomplete information.
\newblock In {\em The 48th Annual Allerton Conference on Communication,
  Control, and Computing (Allerton)}, pages 704--711, June 2010.

\bibitem{yu14a}
H.-F. Yu, C.-J. Hsieh, S.~Si, and I.~S. Dhillon.
\newblock Parallel matrix factorization for recommender systems.
\newblock {\em Knowledge and Information Systems}, 41(3):793--819, 2014.

\bibitem{recht13a}
B.~Recht and C~R{\'e}.
\newblock Parallel stochastic gradient algorithms for large-scale matrix
  completion.
\newblock {\em Mathematical Programming Computation}, 5(2):201--226, 2013.

\bibitem{teflioudi12a}
C.~Teflioudi, F.~Makari, and R.~Gemulla.
\newblock Distributed matrix completion.
\newblock In {\em International Conference on Data Mining (ICDM)}, pages
  655--664, 2012.

\bibitem{lin15a}
A.-Y. Lin and Q.~Ling.
\newblock Decentralized and privacy-preserving low-rank matrix completion.
\newblock {\em Journal of the Operations Research Society of China},
  3(2):189--205, 2015.

\bibitem{boyd06a}
S.~Boyd, A.~Ghosh, B.~Prabhakar, and D.~Shah.
\newblock Randomized gossip algorithms.
\newblock {\em IEEE Transaction on Information Theory}, 52(6):2508--2530, 2006.

\bibitem{mishra14a}
B.~Mishra, G.~Meyer, S.~Bonnabel, and R.~Sepulchre.
\newblock Fixed-rank matrix factorizations and {R}iemannian low-rank
  optimization.
\newblock {\em Computational Statistics}, 29(3--4):591--621, 2014.

\bibitem{dai12a}
W.~Dai, E.~Kerman, and O.~Milenkovic.
\newblock A geometric approach to low-rank matrix completion.
\newblock {\em IEEE Transactions on Information Theory}, 58(1):237--247, 2012.

\bibitem{absil08a}
P.-A. Absil, R.~Mahony, and R.~Sepulchre.
\newblock {\em Optimization Algorithms on Matrix Manifolds}.
\newblock Princeton University Press, Princeton, NJ, 2008.

\bibitem{bonnabel13a}
S.~Bonnabel.
\newblock Stochastic gradient descent on {R}iemannian manifolds.
\newblock {\em IEEE Transactions on Automatic Control}, 58(9):2217--2229, 2013.

\bibitem{huper10a}
K.~H\"uper, U.~Helmke, and S.~Herzberg.
\newblock On the computation of means on {Grassmann} manifolds.
\newblock In {\em International Symposium on Mathematical Theory of Networks
  and Systems (MTNS)}, pages 2439--2441, 2010.

\bibitem{boumal14a}
N.~Boumal, B.~Mishra, P.-A. Absil, and R.~Sepulchre.
\newblock Manopt: a {M}atlab toolbox for optimization on manifolds.
\newblock {\em Journal of Machine Learning Research}, 15(Apr):1455--1459, 2014.

\bibitem{harper15a}
F.~M. Harper and J.~A. Konstan.
\newblock The {MovieLens} datasets: history and contex.
\newblock {\em ACM Transactions on Interactive Intelligent Systems}, 5(4),
  2015.

\end{thebibliography}

\end{document}